\documentclass[11pt]{amsart}
\usepackage{amsopn}                    
\usepackage{amsmath} 
\usepackage{amssymb}
\usepackage[all,2cell,emtex]{xy}
\input xypic                                      
\usepackage[mathscr]{euscript}
\usepackage{mathrsfs}

\theoremstyle{plain} 
\newtheorem{thm}{Theorem}[section]
\newtheorem{prop}[thm]{Proposition}
\newtheorem{lemma}[thm]{Lemma}
\newtheorem{cor}[thm]{Corollary}

\theoremstyle{definition}                                         
\newtheorem{definition}[thm]{Definition}
\newtheorem{paragr}[thm]{}

\newtheorem{rem}[thm]{Remark}

\newtheorem{example}[thm]{Example}

\renewcommand{\mathcal}{\mathscr}
\newcommand{\pref}[1]{{\widehat{ #1 }}}

\newcommand{\To}{\longrightarrow}
\newcommand{\cats}{\varDelta}

\newcommand{\simpl}{\pref{\cats}}

\newcommand{\cat}{{\mathcal{C} \hskip -2.5pt \it{at}}}
\newcommand{\ogrp}{{\omega\text{-}{\mathcal{G} \hskip -2.5pt \it{roupoid}}}}
\newcommand{\set}{{\mathcal{S} \hskip -2.pt \it{et}}}

\newcommand{\Hom}{\mathrm{Hom}}
\newcommand{\sHom}{{\mathcal{H} \hskip -2.5pt \it{om}}}

\newcommand{\nerf}{\mathcal{N}}
\newcommand{\reac}{\mathit{cat}}
\newcommand{\op}[1]{{#1}^{\mathit{op}}}

\newcommand{\sing}{\mathrm{Sing}^{}}

\newcommand{\ho}{\operatorname{\mathsf{Ho}}}

\newcommand{\M}{\mathcal{M}}

\newcommand{\smp}[1]{ \varDelta [#1]}
\newcommand{\cornet}{\varLambda}
\newcommand{\R}{\mathbf{R}}
\newcommand{\derL}{\mathbf{L}}

\newcommand{\s}{\mathbf{s}}

\newcommand{\bord}{\partial}
\newcommand{\Glob}{\mathbb{G}}

\def\TO#1{\mathrel{\hbox to #1pt{\rightarrowfill}}}
\def\OT#1{\mathrel{\hbox to #1pt{\leftarrowfill}}}

\def\limproj{\mathop{\oalign{\rm lim\cr
\hidewidth$\longleftarrow$\hidewidth\cr}}}%
\def\limind{\mathop{\oalign{\rm lim\cr
\hidewidth$\longrightarrow$\hidewidth\cr}}}%
\numberwithin{equation}{thm}

\newcommand{\todouble}{\xymatrixcolsep{1pc}\xymatrix{\ar@<.5ex>[r]\ar@<-.5ex>[r]&}}
\newcommand{\todoubleop}{\xymatrixcolsep{1pc}\xymatrix{\ar@<.5ex>[r]&\ar@<.5ex>[l]}}

\renewcommand{\hookrightarrow}{{\hskip -1.5pt\raise 1.5pt\vbox{\xymatrixcolsep{.9pc}\xymatrix{\ar@{^{(}->}[r]&}}\hskip -3.5pt}}

\newcommand{\intcoin}[5]{\raise 12pt\vbox{\xymatrixcolsep{.9pc}\xymatrixrowsep{.7pc}\xymatrix{%
\scriptstyle #1\ar[r]^{\scriptscriptstyle #5}\ar[d]_{\scriptscriptstyle #4}&\scriptstyle #3\\\scriptstyle #2}}}

\newcommand{\Top}{\mathcal{T}\hskip -2.5pt op}

\newcommand{\ograph}{\omega\text{-}{\mathcal{G} \hskip -2pt \mathit{raph}}}
\newcommand{\ocat}{\omega\text{-}\cat}
\newcommand{\alg}{\mathbf{Alg}}
\newcommand{\linear}[1]{\overline{#1}}
\newcommand{\oper}[1]{{\underline{#1}}}
\newcommand{\joyal}{{\varTheta}}
\newcommand{\height}{\mathrm{ht}}
\newcommand{\tps}{\mathit{top}}

\title[Batanin Higher Groupoids]{Batanin Higher Groupoids\\
and Homotopy Types}
\author{Denis-Charles Cisinski}
\address{Universit\'e Paris 13\\
Institut Galil\'ee\\
LAGA\\
CNRS (UMR 7539)\\
Avenue Jean-Baptiste Cl\'ement\\
93430 Villetaneuse\\
France}
\email{cisinski@math.univ-paris13.fr}
\urladdr{http://www.math.univ-paris13.fr/~cisinski/}
\dedicatory{To Ross Street for his 60th birthday}

\begin{document}
\begin{abstract}
We prove that any homotopy type can be recovered canonically from its
associated weak $\omega$-groupoid. This implies that the homotopy category
of CW-complexes can be embedded in the homotopy
category of Batanin's weak higher groupoids. 
\end{abstract}
\maketitle

\section*{Introduction}

The idea that the homotopy category of CW-complexes
should be equivalent to some reasonnable
homotopy category of $\infty$-groupoids comes from
Grothendieck's letter to Quillen in \emph{Pursuing stacks}~\cite{Gr1}.
But we know that such a result is false as far as we consider
only strict higher groupoids; see e.g. \cite{Ler2,Be2}.
However, we still hope
to prove this equivalence by considering some
weaker notion of higher groupoid.
M.~Batanin~\cite{Ba} did give a more precise
formulation of this: he defined a reasonnable notion
of weak $\omega$-groupoid and a nice functor from
spaces to the category of weak $\omega$-groupoids,
and made the hypothesis that this should give an answer
to this problem; see \cite[Hypothesis page 98]{Ba} for a precise formulation.
Batanin's definition of weak $\omega$-groupoids is based on the notion
of $\omega$-operad. These are operads for which the
operations are not indexed by integers but by (finite planar level) trees.
The underlying algebra has been intensively studied
by R. Street and M. Batanin; see~\cite{Ba,St2,bas}.
C.~Berger~\cite{Be} started to study how one can
define reasonnable homotopy theories from this kind of object.
In particular, thanks to Berger's insights, a good notion
of weak equivalence between weak $\omega$-groupoids
is now available, so that our problem is more precise.
In these notes, we prove that Batanin's functor (or a very
slightly modified version of it)
$$\Pi_\infty : \ho(\Top)\To\ho(\ogrp)$$
is faithful and conservative. This comes from the fact that
one can reconstruct canonically
the homotopy type of a space $X$ from its associated
weak $\omega$-groupoid $\Pi_\infty(X)$ (see \ref{goodgroupoid2}
and its proof). This is proved by relating
Batanin's construction of the functor $\Pi_\infty$
with Berger's results (this is essentially
the contents of Theorem \ref{unicitycanreal} and its corollary).

\section{Omega operads}
\begin{paragr}\label{defomegagraphs}
We first recall the notion of $\omega$-operad (see \cite{Ba,St2,bas}).
Let $\Glob$ be the globe category. The objects of $\Glob$ are the
$\overline{n}$'s for each integer $n\geq 0$.  The globular operators
(i.e. the maps in $\Glob$) are generated by the cosource and cotarget maps
$s,t:\overline{n}\To\overline{n+1}$ subject to the relations $ss=ts$ and $tt=st$.
Recall that an \emph{$\omega$-graph} is a presheaf (of sets) on the category
$\Glob$. The category of $\omega$-graphs is denoted by $\ograph$. The obvious
forgetful functor
$$U:\ocat\To\ograph$$
from the category of \emph{strict} $\omega$-categories to the category
of $\omega$-graphs has a left adjoint
$$F:\ograph\To\ocat \ . $$
The functor $U$ is \emph{monadic}. This means that if $\varepsilon$ denotes the counit of the
adjunction $(F,U)$, the functor
$$Y\longmapsto (UY,U\varepsilon_Y)$$
induces an equivalence of categories $\ocat\simeq\alg_{UF}$.
In other words, if we define $\underline{\omega}=UF$, the category
$\ocat$ is canonically equivalent to the category of algebras on the
monad $\underline{\omega}$ (see for example \cite[Theorem 1.12]{Be}).\\
\indent An \emph{$\omega$-collection} $A$ is a pair $(\underline{A},p_A)$
where $A$ is an endofunctor of the category $\ograph$, and $p_A$ is a
\emph{cartesian} natural transformation
$$p_A : \underline{A}\To\underline{\omega} \ . $$
This means that for any map $X\To Y$ of $\omega$-graphs, the
commutative square
$$\xymatrix{
\underline{A}(X)\ar[r]\ar[d]&\underline{\omega}(X)\ar[d]\\
\underline{A}(Y)\ar[r]&\underline{\omega}(Y)}$$
is a pullback square in $\ograph$. A map of $\omega$-collections
$A\To B$ is a natural transformation $q:\underline{A}\To\underline{B}$
such that $p_B q= p_A$. The category of $\omega$-collections is
endowed with a monoidal structure induced by the composition
of endofunctors: $A\otimes B=(\underline{A}\circ\underline{B},p_{A\otimes B})$,
where $p_{A\otimes B}$ is the composition below.
$$\xymatrix{\underline{A}\circ\underline{B}\ar[r]^{\underline{A}\star p^{}_B}&
\underline{A}\circ\underline{\omega}\ar[r]^{p^{}_A\star\underline{\omega}}&
\underline{\omega}\circ\underline{\omega}\ar[r]^\mu&\underline{\omega}}$$
This tensor product is well defined as the natural transformation $\mu$
above is cartesian (see \cite{St2}).
An \emph{$\omega$-operad} is a monoid object in the category
of $\omega$-collections. In other words, an $\omega$-operad is
an $\omega$-collection $A=(\underline{A},p_A)$ endowed with
a structure of a monad on $\underline{A}$ over $\underline{\omega}$.
For example, $\omega=(\oper{\omega},1_{\oper{\omega}})$
is an $\omega$-operad. This is clearly the \emph{terminal $\omega$-operad},
that is the terminal object in the category of $\omega$-operads.
\end{paragr}

\begin{paragr}\label{deftrees}
Let $X$ be an $\omega$-graph. For $n\geq 0$, we denote
by $X_n$ the set $X(\overline{n})$ of $n$-cells in $X$. 
We define $\leq_X$ to be the preorder
relation on the set $\amalg_{n\geq 0}X_n$
generated by
$$s^*(x)\leq_X x\quad\text{and}\quad x\leq_X t^*(x)\quad\text{for any $x$ in $X_n$, $n\geq 1$,}$$
where $s^*,t^*:X_{n}\To X_{n-1}$ are the maps induced by the operators
$s,t:\overline{n-1}\To\overline{n}$.
We define a (finite planar level) \emph{tree} as a
non empty $\omega$-graph
$T$ such that the set $\amalg_{n\geq 0}T_n$ is finite
and $\leq_T$ is a total order\footnote{This definition of trees
is not the usual one (see for example \cite{Ba,Be}), but
is equivalent to it: finite planar level trees in the usual sense
correspond bijectively to $\omega$-graphs satisfying the
properties listed above (see \cite{St2}).
The definition we have chosen is related with
the canonical embedding of trees in $\omega$-graphs
(see also \cite[p. 61 and p. 80]{Ba}).}.\\
\indent For example, for any $n\geq 0$, the $\omega$-graph
$\overline{n}$ (that is the $\omega$-graph represented
by the corresponding object in $\Glob$) is a tree. We will say
that a tree $T$ is \emph{linear} if $T\simeq\overline{n}$
for some integer $n\geq 0$. A tree is \emph{non linear}
if it is not linear.\\
\indent A tree $T$ is of height $\leq n$ if for any integer $k > n$,
$T^k=\varnothing$. We define the height of $T$ as the integer
$$\height(T)=\mathrm{min}\{n\, | \, \text{$T$ is of height $\leq n$}\} \ . $$

For an $\omega$-graph $X$ and a tree $T$, one defines $X^T$
as the set of maps
$$X^T=\Hom_{\ograph}(T,X) \ . $$

A morphism of trees is just
a morphism of $\omega$-graphs. Any morphism
of trees is a monomorphisms (see for example \cite[Lemma 1.3]{Be}).
We denote by $\joyal_0$ the full subcategory of $\ograph$
whose objects are trees. We will consider $\joyal_0$ as a site with the Grothendieck
topology defined by epimorphic families of maps (see \cite[Definition 1.5]{Be}).
One checks easily that a family of maps $\{T_{i}\To T\}_{i\in I}$ in $\joyal_0$
is a covering if and only if the induced map $\amalg_{i\in I}T_{i}\To T$ is an
epimorphism in the category $\ograph$.
\end{paragr}

\begin{paragr}\label{defoperAnerve}
Let $A$ be an $\omega$-operad. We denote by $\alg_\oper{A}$
the category of algebras over the monad $\oper{A}$. We then have
a free $\oper{A}$-algebra functor
$$\ograph\To\alg_\oper{A}\quad , \qquad X\longmapsto\oper{A}(X)$$
whose right adjoint is the forgetful functor
$$\oper{A}\To\ograph\quad , \qquad C\longmapsto C \ . $$
We define the category $\joyal_A$ as the full subcategory
of $\alg_\oper{A}$ spanned by objects of the form $\oper{A}(T)$
for any tree $T$ (this is the full subcategory of the
Kleisli category of the monad $\oper{A}$
spanned by trees). One has by definition a functor
$$i : \joyal_0\To\joyal_A\quad , \qquad T\longmapsto \oper{A}(T) \ . $$
An \emph{$A$-cellular set} is a presheaf of sets on the category
$\joyal_A$.\\
\indent The inclusion functor of $\joyal_A$ in $\alg_\oper{A}$
induces a nerve functor from $\alg_\oper{A}$ to the category
$\smash{\pref{\joyal}_A=\set^{\op{\joyal}_A}}$ of presheaves of sets on the
category $\joyal_A$
$$\nerf_A : \alg_\oper{A}\To\pref{\joyal}_A \ ,
\quad C\longmapsto \big(\oper{A}(T)\longmapsto C^T=\Hom_{\alg_\oper{A}}(\oper{A}(T),C)\big) \ . $$
The nerve functor $\nerf_A$ is a right adjoint to the left Kan extension
of the inclusion functor of $\joyal_A$ in $\alg_\oper{A}$
$$\reac_A : \pref{\joyal}_A\To\alg_\oper{A} \ . $$
We have finally an inverse image functor
$$i^* : \pref{\joyal}_A\To\pref{\joyal}_0 \ . $$
Following  \cite{Be}, we say that a presheaf $X$ on $\joyal_A$ is a
\emph{$\joyal_A$-model} if its restriction $i^*(X)$ on $\joyal_0$
is a sheaf for the Grothendieck topology defined in \ref{deftrees}.
\end{paragr}

\begin{thm}[C. Berger]\label{opernerve}
The nerve functor $\nerf_A : \alg_\oper{A}\To\pref{\joyal}_A$
is fully faithful. Its essential image consists of the $\joyal_A$-models.
\end{thm}

\begin{proof}
See \cite[Theorem 1.17]{Be}.
\end{proof}

\begin{example}\label{exinitomega}
If $\varnothing$ denotes the initial $\omega$-operad
(that is the corresponding monad is the identity
of $\ograph$), then $\joyal_\varnothing=\joyal_0$ is the already defined
category of trees. By definition, the corresponding category of algebras
is the category of $\omega$-graphs. Hence Theorem \ref{opernerve}
says that the category $\ograph$ is canonically equivalent to the
category of sheaves on $\joyal_0$ for the Grothendieck topology
defined in \ref{deftrees}.
\end{example}

\begin{example}\label{exomega}
If $\omega$ is the terminal $\omega$-operad, the
category $\joyal_\omega$ is the opposite of
Joyal's category of $\omega$-disks (see \cite[Proposition 2.2]{Be}).
By definition, the corresponding category
of algebras is the category of strict $\omega$-categories.
\end{example}

\begin{lemma}\label{Segalcondition0}
A presheaf $X$ on $\joyal_0$ is a sheaf if and only if
for any tree $T$, the map
$$\Hom_{\pref{\joyal}_0}(T,X)\To\limproj_{\linear{n}\To T}\Hom_{\pref{\joyal}_0}(\linear{n},X)$$
is a bijection.
\end{lemma}

\begin{proof}
This comes from the fact that for any tree $T$, the map
$$\limind_{\linear{n}\To T}\linear{n}\To T$$
is an isomorphism in the category of $\omega$-graphs.
\end{proof}

\begin{paragr}\label{defnotareperesacell}
For a tree $T$, we denote by $\joyal_A [T]$ the
presheaf on $\joyal_A$ represented by $\oper{A}(T)$.
Hence for any tree $S$, we have
$$\joyal_A [T] (S)\simeq\Hom_{\joyal_A}(\oper{A}(S),\oper{A}(T)) \ . $$
We also have the identification $\nerf_A(\oper{A}(T))=\joyal_A [T]$.
\end{paragr}

\begin{prop}\label{Segalcondition}
Let $A$ be an $\omega$-operad. Then an $A$-cellular set $X$ is
(isomorphic to) the nerve of an $\oper{A}$-algebra if and only if for any tree $T$, the map
$$\Hom_{\pref{\joyal}_A}(\joyal_A [T],X)\To\limproj_{\linear{n}\To T}\Hom_{\pref{\joyal}_A}(\joyal_A[\linear{n}],X)$$
is a bijection.
\end{prop}

\begin{proof}
If $i_! :\pref{\joyal}_0\To\pref{\joyal}_A$
denotes the left Kan extension of the canonical functor $i$ from $\joyal_0$ to $\joyal_A$,
one has $i_! T\simeq\oper{A}(T)$ for any tree $T$.
This implies that the given condition is equivalent to saying that for any tree $T$, the map
$$\Hom_{\pref{\joyal}_0}(T,i^*(X))\To\limproj_{\linear{n}\To T}\Hom_{\pref{\joyal}_0}(\linear{n},i^*(X))$$
is a bijection. Moreover, we know that $X$ is the nerve of an $\oper{A}$-algebra
if and only if $i^*(X)$ is a sheaf on $\joyal_0$. The result follows by
Lemma \ref{Segalcondition0}.
\end{proof}

\begin{paragr}\label{defbordtree}
For a tree $T$ we define its \emph{boundary} to be the presheaf
on $\joyal_0$
$$\bord T = \bigcup_{S\subset T,\\ S\neq T} S \ , $$
where the $S$'s run over the proper subtrees of $T$.
For an $\omega$-operad $A$, we define $\bord\joyal_A [T]=i_!(\bord T)$.
We then have a canonical inclusion of $A$-cellular sets
$$\bord\joyal_A [T]\To\joyal_A [T] \ . $$
\end{paragr}

\begin{prop}\label{algebrasandboundaries}
Let $A$ be an $\omega$-operad. Then for any $\joyal_A$-model $X$,
and any non linear tree $T$, the canonical map
$$\Hom_{\pref{\joyal}_A}(\joyal_A [T],X)\To
\Hom_{\pref{\joyal}_A}(\bord\joyal_A [T],X)$$
is an isomorphism.
\end{prop}

\begin{proof}
It is sufficient to prove that the canonical map
$$\limind_{S\subset T , \\ S\neq T}\oper{A}(S)
\simeq\oper{A}\big(\limind_{S\subset T , \\ S\neq T}S\big)\To \oper{A}(T)$$
is an isomorphism in the category of $\oper{A}$-algebras.
As the functor $\oper{A}$ from $\omega$-graphs to $\oper{A}$-algebras
preserves colimits, we can suppose for this that $A$ is the initial $\omega$-operad.
But the inclusions $S\subset T$, $S\neq T$, form an epimorphic family
in the category of presheaves over $\joyal_0$, which implies that $\bord T$
is a covering sieve of $T$. This proves that the sheaf associated to the presheaf $\bord T$
is canonically isomorphic to $T$. Hence the result.
\end{proof}

\section{Contractible $B$-operads}

\begin{paragr}\label{batreal}
Let $\Top$ be the category of (compactly generated) topological
spaces. We define a functor
\begin{equation}\label{globtopolograph}
b : \Glob \To\Top
\end{equation}
by $b(\linear{n})=B^n_\tps$, where $B^n_\tps$ denotes the
euclidian $n$-dimensional ball. The operator $s$ (resp. $t$)
from $\linear{n-1}$ to $\linear{n}$ is sent to the continuous
map $s_\tps$ (resp. $t_\tps$) which sends $B^{n-1}_\tps$
to the north hemisphere of $B^n_\tps$ (resp. to the south hemisphere
of $B^n_\tps$); see \cite[Example 2, p. 58]{Ba}.
This defines a functor
\begin{equation}\label{globtopolograph2}
b^*:\Top\To\ograph\quad , \qquad X\longmapsto (\linear{n}\longmapsto\Hom_{\Top}(B^n_\tps,X)) \ . 
\end{equation}
The aim of this section is to explain how we can produce
some $\omega$-operads $A$ such that for any topological space $X$,
the $\omega$-graph $b^*(X)$ is endowed with a structure of an
$\oper{A}$-algebra\footnote{This construction is due to Batanin~\cite{Ba}.}.
Taking the left Kan extension of $b$ leads to a left adjoint to $b^*$,
that is the unique cocontinuous functor
$$b_! : \ograph\To\Top$$
whose restriction to $\Glob$ is $b$. By restriction to $\joyal_0$, this defines a functor
$$R_0 :\joyal_0\To\Top\quad , \qquad T\longmapsto B^T_\tps=b_!(T) \  . $$
More explicitely, one has a canonical isomorphism
$$\limind_{\linear{n}\To T} B^n_\tps\simeq B^T_\tps \ . $$
\end{paragr}

\begin{paragr}\label{defbataninoperad}
According to \cite[Proposition 7.2]{Ba}, the coglobular topological
space $b$ defines an $\omega$-operad $B=\op{E}(b)$.
We will call $B$ the \emph{topological $\omega$-operad}
(see also \cite[section 9, p. 97 sq]{Ba}). The corresponding
monad $\oper{B}$ on $\ograph$ is defined as follows.
For an $\omega$-graph $X$ and an integer $n$,
the $\omega$-graph $\oper{B}(X)$ is given by the formula
\begin{equation}\label{deffreeBalg}
\oper{B}(X)_n=\coprod_{\height(T)\leq n}\Hom_{}(B^n_\tps,B^T_\tps)\times X^T
\end{equation}
where $\Hom_{}(B^n_\tps , B^T_\tps)$ denotes the set of maps
between the corresponding coglobular spaces\footnote{Any tree has a canonical
structure of a coglobular $\omega$-graph (see e.g. \cite[Definition 1.8]{Be}).
Hence, by functoriality, for any tree $T$, $B^T_\tps$
is canonically endowed with a structure of a coglobular space.}.
The natural transformation
$$p_B : \oper{B}\To\oper{\omega}$$
is defined by the projections
$$\Hom_{}(B^n_\tps,B^T_\tps)\times X^T\To X^T$$
once we remember the canonical identification
(see \cite[Definition 1.8 and Theorem 1.12]{Be})
\begin{equation}\label{defreefomega}
\oper{\omega}(X)_n\simeq\coprod_{\height(T)\leq n} X^T \ . 
\end{equation}
\end{paragr}

\begin{prop}\label{realBalg}
There is a canonical functor
$$R_B : \joyal_{B}\To\Top$$
such that the following diagram commutes.
$$\xymatrix{
\joyal_{0}\ar[r]^{R_0}\ar[d]_{i}&\Top\\
\joyal_{B}\ar[ur]_{R_B}
}$$
\end{prop}

\begin{proof}
As the functor $i$ is bijective on objects, we only have to define
the functor $R_B$ on morphisms.
Let $n\geq 0$ be an integer, and let $T$ be a tree.
We have
$$\begin{aligned}
\Hom_{\joyal_B}(\oper{B}(\linear{n}),\oper{B}(T))
&\simeq \oper{B}(T)_n\\
&\simeq\coprod_{\height(U)\leq n}
\Hom_{}(B^n_\tps,B^U_\tps)\times \Hom_{\joyal_0}(U,T)
\end{aligned}$$
hence the maps
$$\Hom_{}(B^n_\tps,B^U_\tps)\times \Hom_{\joyal_0}(U,T)
\To\Hom_{\Top}(B^n_\tps,B^T_\tps)$$
defined by $(g,f)\longmapsto R_0(f)\circ g$
induce a map
$$\Hom_{\joyal_B}(\oper{B}(\linear{n}),\oper{B}(T))
\To\Hom_{\Top}(B^n_\tps,B^T_\tps) \ .$$
For two trees $S$ and $T$, we obtain a map
$$\limproj_{\linear{n}\To S}\Hom_{\joyal_B}(\oper{B}(\linear{n}),\oper{B}(T))
\To\limproj_{\linear{n}\To S}\Hom_{\Top}(B^n_\tps,B^T_\tps)$$
that is a map
$$\Hom_{\joyal_B}(\oper{B}(S),\oper{B}(T))
\To\Hom_{\Top}(B^S_\tps,B^T_\tps) \ . $$
This defines the functor $R_B$ and ends the proof.
\end{proof}

\begin{paragr}\label{defBcellulartopreal}
We define the functor
$$|~-~|_B : \pref{\joyal}_B\To\Top\quad , \qquad
X\longmapsto |X|_B$$
as the left Kan extension of the functor
$R_B$ given by Proposition \ref{realBalg}. It has a right adjoint
$$\Pi^B_\infty : \Top\To\pref{\joyal}_B\quad ,
\qquad X\longmapsto \big(\oper{B}(T)\longmapsto\Hom_{\Top}(B^T_\tps,X) \big) \ . $$
\end{paragr}

\begin{prop}\label{piinftyBalg}
For any (compactly generated) topological space $X$,
the $B$-cellular set $\Pi^B_\infty(X)$ is a $\oper{B}$-algebra.
\end{prop}

\begin{proof}
According to Proposition \ref{Segalcondition}, it is sufficient to check
that for any tree $T$, the map
$$\Hom_{\pref{\joyal}_B}(\joyal_B [T],\Pi^B_\infty(X))\To
\limproj_{\linear{n}\To T}
\Hom_{\pref{\joyal}_B}(\joyal_B [\linear{n}],\Pi^B_\infty(X))$$
is a bijection. The canonical identifications
$$\begin{aligned}
\limproj_{\linear{n}\To T}
\Hom_{\pref{\joyal}_B}(\joyal_B [\linear{n}],\Pi^B_\infty(X))
&\simeq\limproj_{\linear{n}\To T}\Hom_{\Top}(B^n_\tps,X)\\
&\simeq\Hom_{\Top}(B^T_\tps,X)\\
&\simeq\Hom_{\pref{\joyal}_B}(\joyal_B [T],\Pi^B_\infty(X))
\end{aligned}$$
end this proof.
\end{proof}

\begin{cor}\label{topbordball}
For any non linear tree $T$, the continuous map
$$|\bord\joyal_B [T]|_B\To |\joyal_B [T]|_B=B^T_\tps$$
is an homeomorphism.
\end{cor}

\begin{proof}
By Theorem \ref{opernerve} and Proposition \ref{Segalcondition},
this is reformulation of Propositions \ref{algebrasandboundaries} and \ref{piinftyBalg}
and of the Yoneda Lemma applied to $\Top$.
\end{proof}

\begin{definition}\label{defgentopoperad}
A \emph{$B$-operad} is an $\omega$-operad $A$
endowed with a morphism of $\omega$-operads $\varphi^{}_A : A\To B$
from $A$ to the topological $\omega$-operad $B$.
\end{definition}

\begin{rem}\label{defBrealBoper}
For a $B$-operad $A$, the morphism $\varphi^{}_A$ defines a functor
$$\varphi^{}_A : \joyal_A\To\joyal_B\quad , \qquad \oper{A}(T)\longmapsto\oper{B}(T) \ . $$
We define the functor
\begin{equation}\label{defRA}
R_A : \joyal_A\To\Top
\end{equation}
as the composition of $\varphi^{}_A$ and of the functor $R_B$ (see \ref{defBcellulartopreal}).
In other words, for any tree $T$, we have $R_A(\oper{A}(T))=R_B(\oper{B}(T))=B^T_\tps$.
We define the functor
\begin{equation}\label{defrealRA}
|~-~|_A : \pref{\joyal}_A\To\Top
\end{equation}
as the left Kan extension of $R_A$. It has a right adjoint
\begin{equation}\label{defPiA0}
\Pi^A_\infty : \Top\To\pref{\joyal}_A
\end{equation}
defined by $\Pi^A_\infty=\varphi^*_A\Pi^B_\infty$, where
$\varphi^*_A : \pref{\joyal}_B\To\pref{\joyal}_A$
is the inverse image functor of the functor $\varphi^{}_A$ above.
As the functor $\varphi^*_A$ obviously sends
$\joyal_B$-models on $\joyal_A$-models, it follows from Theorem \ref{opernerve}
and Proposition \ref{piinftyBalg} that $\Pi^A_\infty$ can be defined as a functor
\begin{equation}\label{defPiA}
\Pi^A_\infty : \Top\To\alg_{\oper{A}} \ .
\end{equation}
We remark that for any topological space $X$, the
underlying $\omega$-graph of $\Pi^A_\infty(X)$
is the $\omega$-graph $b^*(X)$ introduced in \ref{batreal}.
\end{rem}

\begin{paragr}\label{defcontractoperad}
Following~\cite[Definition 1.20]{Be}\footnote{The definition
of contractible $\omega$-operads in terms of lifting properties
in fact appeared first in Lemma 1.1 from \cite{Ba2}.},
we say that an $\omega$-operad $A$ is \emph{contractible}
if for any $\omega$-graph $X$, the canonical map
$\oper{A}(X)\To\oper{\omega}(X)$ has the right lifting property
with respect to inclusions $\bord\linear{n}\To\linear{n}$, $n\geq 0$.
In other words, we ask that any solid commutative square
\begin{equation}\label{contractlift}
\begin{split}\xymatrix{
\bord\linear{n}\ar[r]\ar[d]&\oper{A}(X)\ar[d]\\
\linear{n}\ar[r]\ar@{-->}[ur]^f&\oper{\omega}(X)}
\end{split}\end{equation}
admits a diagonal filler $f$ such that the diagram still commutes.
As the natural transformation $\oper{A}\To\oper{\omega}$
is cartesian, it is sufficient to check
the lifting property above in the case where $X=e$ is the terminal
$\omega$-graph (see \cite[Definition 1.20]{Be}).
It is then easy to check that this definition
is somehow equivalent to the notion of \emph{$\omega$-operad with
contractions} as defined in
\cite{Leibook,Lei}\footnote{However it is not obvious
at all that this definition is equivalent to Batanin's definition of
contractible $\omega$-operad with system of
compositions~\cite[Definitions 8.2 and 8.4]{Ba}.
Any contractible $\omega$-operad as defined in \ref{defcontractoperad}
is a contractible $\omega$-operad with system of compositions in Batanin's sense,
but the converse is still a conjecture.}. More precisely, an $\omega$-operad
with contractions as defined by Leinster corresponds to a contractible $\omega$-operad $A$
with a given lift $f$ for all the commutative squares of type \eqref{contractlift}
where $X$ is the terminal $\omega$-graph.
\end{paragr}

\begin{prop}\label{battopcontr}
The topological $\omega$-operad $B$ is contractible.
\end{prop}

\begin{proof}
Let $X$ be an $\omega$-graph. For any integer $n\geq 0$, we have by construction
$$\Hom_{\ograph}(\linear{n},\oper{B}(X))\simeq
\oper{B}(X)_n=\coprod_{\height(T)\leq n}\Hom_{}(B^n_\tps,B^T_\tps)
\times X^T  \ . $$
We deduce from this that
$$\Hom_{\ograph}(\bord\linear{n},\oper{B}(X))\simeq
\coprod_{\height(T)\leq n}\Hom_{}(\bord B^n_\tps,B^T_\tps)\times X^T$$
where $\bord B^n_\tps=S^{n-1}$ is the $n-1$-dimensional sphere.
Using the formula \eqref{defreefomega}, one concludes that it is sufficient
to check that any map from $S^{n-1}$ to $B^T_\tps$
can be extended to $B^n_\tps$. This follows from the fact that $B^T_\tps$
is contractible.
\end{proof}

\begin{thm}[T. Leinster]\label{existinitcontractoper}
There is a universal\footnote{This is in the sense defined
in \cite[Definition 8.5]{Ba}: universal means weakly initial.}
contractible $\omega$-operad $K$.
\end{thm}

\begin{proof}
By \cite[Proposition 9.2.2 and Appendix G]{Leibook},
the category of $\omega$-operads with contractions
as defined in \cite[Definition 9.2.1]{Leibook} has
an initial object $K$. The assertion now comes from the fact that
any contractible $\omega$-operad $A$ can be endowed
with the structure of an $\omega$-operad with contractions
(just chose the lift $f$ for all the commutative squares of type \eqref{contractlift}
where $X$ is the terminal $\omega$-graph),
so that $K$ is also a weak initial object in the category of
contractible $\omega$-operad.
\end{proof}

\begin{rem}
This theorem has to be compared with Batanin's result \cite[Theorem 8.1]{Ba}
saying that there is a universal contractible $\omega$-operad with system of compositions.
\end{rem}

\begin{paragr}\label{defweakcat}
It is then reasonnable to define a \emph{weak $\omega$-category}
as a $\oper{K}$-algebra (and any such weak $\omega$-category
is a weak $\omega$-category in Batanin's sense~\cite{Ba},
but we don't know if the converse is true).\\
\indent As the topological $\omega$-operad $B$ is contractible, one
has a map of $\omega$-operads $\varphi^{}_K : K\To B$.
In other words, $K$ is a $B$-operad (see \ref{defgentopoperad}).
It is obvious that $K$ is also the initial contractible $B$-operad.
As explained in Remark \ref{defBrealBoper},
the map $\varphi^{}_K$ induces a canonical functor
\begin{equation}\label{defKtoB}
\varphi^{}_K : \joyal_K\To\joyal_B \ . 
\end{equation}
as well as a functor
\begin{equation}\label{deffondinftygroupoid}
\Pi^K_\infty : \Top\To\alg_\oper{K} \ .
\end{equation}
\end{paragr}

\section{Topological realizations of cellular spaces}

\begin{paragr}
Let $\cats$ be the category of simplices.
The objects are the totally ordered sets
$[n]=\{ 0, \ldots , n \}$
for any integer $n\geq 0$, and the morphisms are the
order-preserving maps. The category $\smash{\simpl}$ of simplicial sets
is the category of presheaves of sets on $\cats$.
The simplicial set represented by $[n]$
is denoted by $\smp{n}$, and the boundary of
$\smp{n}$ is denoted by $\bord\smp{n}$. For
any integers $0\leq k\leq n$, $n\geq 1$, $\cornet^k[n]$
is the $k$-horn of  $\smp{n}$ (see \cite{GZ,GJ2}).
\end{paragr}

\begin{paragr}\label{defnotagencofib}
Let $A$ be an $\omega$-operad.
We denote by $\s\pref{\joyal}_A$
the category of simplicial $A$-cellular sets.
We first recall here how one can produce
a reasonnable homotopy theory of
simplicial $A$-cellular sets following
C.~Berger's construction~\cite{Be}.
We define three sets of maps of simplicial
$A$-cellular sets as follows.
The set $I^{}_A$ consists of the inclusions
$$\bord\joyal_A [T]\otimes\smp{n}
\cup\joyal_A [T]\otimes\bord\smp{n}\To
\joyal_A [T]\otimes\smp{n}$$
for any tree $T$ and any integer $n\geq 0$.
The set $J^\prime_A$ consists of the inclusions
$$\bord\joyal_A [T]\otimes\smp{n}
\cup\joyal_A [T]\otimes\cornet^k[n]\To
\joyal_A [T]\otimes\smp{n}$$
for any tree $T$ and any integers $0\leq k\leq n$, $n\geq 1$.
Finally, the set $J^{\prime\prime}_A$ is consists of the inclusions
$$\joyal_A[S]\otimes\smp{n}\cup\joyal_A[T]\otimes\bord\smp{n}
\To\joyal_A[T]\otimes\smp{n}$$
for any map of trees $S\To T$ and any integer $n\geq 0$.
\end{paragr}

\begin{thm}[C. Berger]\label{thmberger1}
Let $A$ be an $\omega$-operad. The category
of simplicial $A$-cellular sets is endowed with
a left proper cofibrantly generated model category
structure. The cofibrations are generated by the set
$I^{}_A$, and the trivial cofibrations are generated
by the set of maps $J^\prime_A\cup J^{\prime\prime}_A$.
\end{thm}

\begin{proof}
See \cite[Proposition 4.11]{Be}.
\end{proof}

\begin{rem}\label{remthmeberger1}
We did not say how the weak equivalences
are defined: this is not necessary as the cofibrations
and the trivial cofibrations of this model structure are
defined so that all the structure is already determined;
see \cite{Be}. However, a way to define this model structure
is to prove that the category $\smash{\s\pref{\joyal}_A}$
is endowed with a proper cofibrantly generated
model category structure with the termwise
simplicial weak equivalences as weak equivalences.
For this model structure, the cofibrations
are generated by $I^{}_A$, and the trivial cofibrations
are generated by $J^\prime_A$. The model structure
of Theorem \ref{thmberger1} can be defined as the left Bousfield
localization of the latter by the maps $S\to T$ for any trees $S$ ant $T$
(this is the precise formulation given in \cite[Proposition 4.11]{Be}).
\end{rem}

\begin{definition}\label{defcanrealoper}
Let $A$ be an $\omega$-operad. A \emph{canonical
realization of $A$} is a functor
$$R : \joyal_A\To \Top$$
satisfying the following properties.
\begin{itemize}
\item[R1] For any tree $T$, the space $R(\oper{A}(T))$
is (weakly) contractible.
\item[R2] If $R_!=|~-~|_R : \pref{\joyal}_A\To\Top$ denotes the left Kan extension
of $R$, then for any tree $T$, the canonical continuous map
$$|\bord\joyal_A [T]|_R\To |\joyal_A [T]|_R=R(\oper{A}(T))$$
is a cofibration of topological spaces (for the usual
model category structure on $\Top$).
\end{itemize}
\end{definition}

\begin{example}\label{bataninrealcan}
For any $B$-operad $A$, the functor $R_A$~\eqref{defRA}
is a canonical realization of $A$: the spaces $B^T_\tps$
are contractible CW-complexes, and for any tree $T$,
the map
$$|\bord\joyal_A [T]|_A\To |\joyal_A [T]|_A=B^T_\tps$$
is a cofibration. To see this, as any homeomorphism
is a cofibration of topological spaces, it is sufficient
to check this property when $T$ is linear (Corollary \ref{topbordball}).
But if $T=\linear{n}$, then $|\bord\joyal_A [T]|_A$
is the $n-1$-dimensionnal sphere $S^{n-1}$, and the
canonical inclusion of $S^{n-1}$ in the $n$-dimensionnal ball
is obviously a cofibration.
\end{example}

\begin{example}\label{exdisks}
As the category $\joyal_\omega$ is the opposite category
of Joyal's category of $\omega$-disks,
the canonical $\omega$-disk structure on the coglobular
space $b$~\eqref{globtopolograph} defines a functor
$$D_\omega : \joyal_\omega\To\Top$$
that happends to be a canonical realization of the
terminal $\omega$-operad (see \cite[Propositions 2.2 and 2.6]{Be}).\\
\indent If $A$ is an $\omega$-operad, we have a canonical functor
$$p : \joyal_A\To\joyal_\omega\quad , \qquad
\oper{A}(T)\To\oper{\omega}(T) \ . $$
We can define a functor
$$D_A : \joyal_A\To \Top$$
by the formula $D_A = D_\omega\circ p$.
One checks immediately that $D_A$ is a canonical realization
of $A$. We will call $D_A$ the \emph{$\omega$-disks
realization} of $A$.
\end{example}

\begin{example}\label{cancanreal}
Let $A$ be an $\omega$-operad.
We define the \emph{categorical realization} $H$ of $A$
as follows. For this, we recall the following general construction.
Given a (small) category $C$ and a presheaf of sets $X$
on $C$, we define the \emph{category of elements of $X$}
as the category $C/X$ whose objects are pairs
$(c,s)$, where $c$ is an object of $C$, and
$s$ a section of $X$ over $c$. A map in $C/X$
from $(c,s)$ to $(c',s')$ is a map $u:c\To c'$ in $C$
such that $u^*(s')=s$. We have a canonical forgetful functor
$$C/X\To C\quad , \qquad (c,s)\longmapsto c \ . $$
\indent For a cellular space $X$, we define the
category $\joyal_A /X$ as the category of elements of $X$.
This defines a functor
from the category of cellular sets to the category
of small categories. Taking the topological realization of
the nerve of the $\joyal_A / X$'s thus defines a functor
$$|~-~|_{H} : \pref{\joyal}_A\To\Top \ . $$
One can show that this functor is cocontinuous
and sends monomorphisms of $A$-cellular sets
to cofibrations of topological spaces: it follows from
\cite[Corollaire A.1.12]{Ci3} that the functor
$$X\longmapsto\text{simplicial nerve of $\joyal_A /X$}$$
is cocontinuous and preserves monomorphisms,
so that this is a consequence of the fact that
the topological realization functor of simplicial sets
is cocontinuous and sends monomorphisms to
cofibrations of topological spaces.
We conclude that the functor
$$H : \joyal_A\To\Top\quad , \qquad
\oper{A}(T)\longmapsto |\joyal_A[T]|_{H}$$
is a canonical realization of $A$: we have already verified R2,
and for R1, it is obvious that $|\joyal_A[T]|_H$ is contractible
as it is the realization of a category with a terminal object.
\end{example}

\begin{paragr}\label{canrealsimpl}
Given a canonical realization $R$ of an $\omega$-operad $A$,
we get a functor
$$\joyal_A\times\cats\To\Top\quad , \qquad
(\oper{A}(T) , [n])\longmapsto R(\oper{A}(T))\times\varDelta^n_\tps$$
where $\varDelta^n_\tps$ is the topological $n$-simplex.
The left Kan extension of the latter is a cocontinuous functor
$$|~-~|^\s_R : \s\pref{\joyal}_A\To\Top$$
whose right adjoint
$$\sing_R : \Top\To\s\pref{\joyal}_A$$
is defined by
$$X\longmapsto(\oper{A}(T) , [n])\longmapsto
\Hom_{\Top}(R(\oper{A}(T))\times\varDelta^n_\tps,X)) \ . $$
\end{paragr}

\begin{prop}\label{quillencanreal}
For any canonical realization $R$ of an $\omega$-operad $A$,
the functor
$$|~-~|^\s_R : \s\pref{\joyal}_A\To\Top$$
is a left Quillen functor.
\end{prop}

\begin{proof}
It is sufficient to check that the functor $| - |_R$
sends the elements of $I^{}_A$
(resp. of $J^\prime_A\cup J^{\prime\prime}_A$)
to cofibrations (resp. to trivial cofibrations) of topological
spaces. As the model category of topological spaces
is monoidal with respect to the cartesian product, it is
sufficient to check the following facts.
\begin{itemize}
\item[(a)] For any integer $n\geq 0$, the topological realization
of the inclusion $\bord\smp{n}\To\smp{n}$ is a cofibration of
topological spaces.
\item[(b)] For any tree $T$, the map
$|\bord\joyal_A[T]|_R\To|\joyal_A[T]|_R$ is a cofibration of
topological spaces.
\item[(c)] For any integers $0\leq k\leq n$, $n\geq 1$,
the topological realization of $\cornet^k[n]$ is contractible.
\item[(d)] For any map of trees $S\to T$, the map
$|\joyal_A[S]|_R\To|\joyal_A[T]|_R$ is a weak equivalence.
\end{itemize}
Properties (a) and (c) are well known. Property (b)
(resp. (d)) follows from con\-di\-tion R2 (resp. R1) of
Definition \ref{defcanrealoper}.
\end{proof}

\begin{paragr}\label{defnotaderivreal}
We denote by
$$|~-~|^\derL_R : \ho(\s\pref{\joyal}_A)\To\ho(\Top)$$
the total left derived functor of the functor $|~-~|^\s_R$.
\end{paragr}

\begin{paragr}\label{recallhtycolimit}
Let $\M$ be a model category. For a given small category $I$,
let $\M^I$ be the category of functors from $I$
to $\M$. We recall the following facts
about homotopy colimits in model categories
(see e.g. \cite{CS,Ci1,DHK}). The colimit functor
$$\limind_I : \M^I\To\M$$
has a total left derived functor
$$\derL\limind_I : \ho(\M^I)\To\ho(\M)$$
where $\ho(\M^I)$ denotes the localization of $\M^I$
by the class of termwise weak equivalences. The functor
$\derL\limind_I$ is a left adjoint to the
functor
$$\ho(\M)\To\ho(\M^I)\quad ,
\qquad X\longmapsto X_I$$
(where $X_I$ is the constant diagram with value $X$).
Moreover, any left Quillen functor
$$F : \M\To\M' \ , $$
induces a functor
$$F : \M^I\To{\M}^{\prime \, I}$$
that has a total left derived functor (see \cite[Proposition 6.12]{Ci1})
$$\derL F : \ho(\M^I)\To\ho(\M^{\prime \, I}) \ . $$
By virtue of \cite[Proposition 6.12]{Ci1}, we also
know that if $F$ is cocontinuous,
the functors $\derL F$ preserve homotopy colimits in
the following sense. For any functor $\Phi$ from $I$ to $\M$,
the canonical map
$$\derL\limind_I\derL F(\Phi)\To\derL F(\derL\limind_I \Phi)$$
is an isomorphism in $\ho(\M')$.
\end{paragr}

\begin{paragr}\label{abstractreal}
Let $A$ be an $\omega$-operad. Given a functor
$$U:\joyal_A\To\Top \ , $$
one defines a functor
$$U^\derL : \s\pref{\joyal}_A\To\ho(\Top)$$
as follows. First of all, the functor $U$
induces a functor
$$U^\dagger : \joyal_A\times\cats\To\Top$$
defined by $U^\dagger(\oper{A}(T),[n])=
U(\oper{A}(T))\times\varDelta^n_\tps$.
Given a simplicial $A$-cellular set $X$,
one has its category of elements as a presheaf over
$\joyal_A\times\cats$
$$E(X)=(\joyal_A\times\cats)/X$$
(see example \ref{cancanreal}). We have a canonical functor
$\varphi^{}_X$ from $E(X)$ to $E(e)=\joyal_A\times\cats$ (where $e$
denotes the terminal simplicial $A$-cellular set).
We define the functor
$$U^\dagger_X : E(X)\To\Top$$
as the composition of $\varphi^{}_X$ with $U^\dagger$.
The functor $U^\derL$ is at last defined by the
formula
$$U^\derL(X)=\protect{\derL\limind}_{E(X)}U^\dagger_X \ . $$
It is obvious that if $V$ is another functor from
$\joyal_A$ to $\Top$, any termwise weak equivalence
$U\To V$ induces an isomorphism of functors
$U^\derL\simeq V^\derL$.\\
\indent We define the functor
\begin{equation}\label{equarealabstract}
||~-~|| : \s\pref{\joyal}_A\To\ho(\Top)
\end{equation}
by the formula $||X||=P^\derL(X)$,
where $P:\joyal_A\To\Top$ is the functor
which sends every object of $\joyal_A$ to the
terminal topological space.\\
\indent Let $R$ be a canonical realization of $A$.
One interesting fact is that one has a canonical isomorphism
\begin{equation}\label{equarealabstract2}
R^\derL\simeq ||~-~||
\end{equation}
(thanks to the condition R1 of Definition \ref{defcanrealoper}).
Moreover, we also have the following.
\end{paragr}

\begin{prop}\label{hptycanrealeqabstrreal}
For any simplicial $A$-cellular set $X$, one has a natural isomorphism
$$R^\derL(X)\simeq |X|^\derL_R \ . $$
\end{prop}

\begin{proof}
Recall from \ref{abstractreal} that $E(X)$ is the
category of elements of $X$. By composing the
canonical functor from $E(X)$ to $\joyal_A\times\cats$
with the Yoneda embedding into $\s\pref{\joyal}_A$,
we get a functor
$$h_X : E(X)\To\s\pref{\joyal}_A \ . $$
We have a canonical map
$$z_X : \protect{\derL\limind}_{E(X)}h_X\To\limind_{E(X)}h_X\simeq X$$
in $\ho(\s\pref{\joyal}_A)$. We claim that this map is an isomorphism.
This is sufficient to prove our result: once we know that
$z_X$ is an isomorphism, we have
$$\begin{aligned}
R^\derL(X)&=\protect{\derL\limind}_{E(X)}|h_X|^\derL_R\\
&\simeq|\protect{\derL\limind}_{E(X)}h_X|^\derL_R\quad
\text{(by the end of \ref{recallhtycolimit})}\\
&\simeq|X|^\derL_R \ . 
\end{aligned}$$
Hence we just have to prove our claim that the map
$z_X$ above is an isomorphism. For this, we need some more
abstract results. By \cite[Proposition 4.11]{Be},
the category of simplicial $A$-cellular sets
is endowed with a (left) proper cofibrantly
generated model category category structure
whose weak equivalences are the termwise simplicial
weak equivalences. The cofibrations (resp. the trivial cofibrations)
are generated by the set $I^{}_A$ (resp. by the set $J^\prime_A$);
see \ref{defnotagencofib}.
It follows that the identity of \smash{$\s\pref{\joyal}_A$}
is a left Quillen functor from the latter model structure
to the model structure of Theorem \ref{thmberger1}.
Hence it is sufficent to check that $z_X$ is an isomorphism
in the localization of \smash{$\s\pref{\joyal}_A$} by the
termwise simplicial weak equivalences. But as
the definition of the homotopy category of a model category
and of the homotopy colimit functors \smash{$\derL\limind$}
only depend on the class of weak equivalences, it is sufficient
to prove our claim by considering our favorite
model category structure on \smash{$\s\pref{\joyal}_A$}
for which the class of weak equivalences consists of the
termwise simplicial weak equivalences.
If we take the model structure where the cofibrations
are the monomorphisms, then the fact that $z_X$
is an isomorphism comes from \cite[Proposition 4.4.24]{Ci3}.
(But if the reader prefers the model structure where
the fibrations are the termwise Kan fibrations,
we can also refer to \cite[Proposition 2.9]{dug1}.)
\end{proof}

\begin{thm}\label{unicitycanreal}
Let $A$ be an $\omega$-operad. If $R_1$ and $R_2$
are two canonical realizations of $A$, then the
corresponding total left derived functors
$|~-~|^\derL_{R_1}$ and $|~-~|^\derL_{R_2}$
are isomorphic.
\end{thm}

\begin{proof}
For any simplicial $A$-cellular set $X$, we have
$$|X|^\derL_{R_1}\simeq ||X|| \simeq |X|^\derL_{R_2}$$
by Proposition \ref{hptycanrealeqabstrreal}.
\end{proof}

\begin{cor}\label{controperquillenequ}
Let $A$ be a contractible $\omega$-operad.
Then for any canonical realization $R$ of $A$,
the functor
$$|~-~|^\s_R : \s\pref{\joyal}_A\To\Top$$
is a left Quillen equivalence.
\end{cor}

\begin{proof}
Given a canonical realization $R$ of $A$,
the left Quillen functor $|~-~|^\s_\R$ is a Quillen
equivalence if and only if its total left derived
functor is an equivalence of categories.
Hence, by Theorem \ref{unicitycanreal},
it is sufficient to prove that there exists a canonical realization
$R$ of $A$ such that $|~-~|^\s_R$ is a Quillen equivalence.
But we know this is the case when $R=D_{A}$ is the $\omega$-disks
realization of $A$~(\ref{exdisks}) by virtue of
\cite[Theorems 3.9 and 4.14]{Be}.
\end{proof}

\section{The weak $\omega$-groupoid functor}

\begin{thm}[C. Berger]\label{thmberger2}
Let $A$ be an $\omega$-operad. Then the category
$\s\alg_\oper{A}$ of simplicial $\oper{A}$-algebras
is a left proper cofibrantly generated
model category with the following definition:
a map of simplicial $\oper{A}$-algebras is
a weak equivalence (resp. a fibration) if
its nerve (\ref{defoperAnerve}) is a
weak equivalence (resp. a fibration)
of simplicial $A$-cellular sets. Moreover, the
nerve functor
$$\nerf_\oper{A} : \s\alg_\oper{A}\To \s\pref{\joyal}_A$$
is a right Quillen equivalence.
\end{thm}

\begin{proof}
See \cite[Theorem 4.13]{Be}.
\end{proof}

\begin{paragr}\label{htyAalgebras}
Let $A$ be an $\omega$-operad.
We say that a map of $\oper{A}$-algebras is
a \emph{weak equivalence} if it is a weak equivalence
of simplicial  $\oper{A}$-algebras. We denote by
$\ho(\alg_\oper{A})$ the localization of the category
$\alg_\oper{A}$ by the weak equivalences.
The canonical functor from $\alg_\oper{A}$
to $\s\alg_\oper{A}$ induces a functor
\begin{equation}\label{berger3}
\ho(\alg_\oper{A})\To\ho(\s\alg_\oper{A}) \ .
\end{equation}
It follows from the following result that this latter
functor is an equivalence of categories.
\end{paragr}

\begin{prop}[C. Berger]\label{Acellfirbant}
Let $A$ be an $\omega$-operad.
Then for any fibrant simplicial $A$-cellular set
for the model category structure of Theorem \ref{thmberger1},
the canonical map $X_0\To X$
is a weak equivalence (where $X_0$ denotes
the $A$-cellular set obtained from the evaluation of $X$ at $\smp{0}$,
seen as a discrete simplicial $A$-cellular set in the canonical way).
\end{prop}

\begin{proof}
See \cite[Proposition 4.17]{Be}.
\end{proof}

\begin{thm}\label{goodgroupoid1}
Let $A$ be a contractible $B$-operad.
Then the functor
$$\Pi^A_\infty : \Top\To\alg_\oper{A}$$
preserves weak equivalences, and the
induced functor
$$\Pi^A_\infty : \ho(\Top)\To\ho(\alg_\oper{A})$$
is an equivalence of categories.
\end{thm}

\begin{proof}
For two compactly generated topological spaces
$X$ and $Y$, we denote by $\sHom(X,Y)$
the space of continuous maps
from $X$ to $Y$ endowed with the compact-open topology.
If $X$ is an object of $\Top$, we get a simplicial space
$S(X)$ defined by
$$S(X)_n=\sHom(\varDelta^n_\tps , X) \ . $$
Applying the functor $\Pi^A_\infty$ termwise to this
simplicial space defines a simplicial $\oper{A}$-algebra $\Pi^A_\infty S(X)$.
One can see easily that the functor
$X\longmapsto\nerf_\oper{A}\Pi^A_\infty S(X)$
is a right adjoint to the left Quillen functor $|~-~|^\s_R$
where $R=R_A$ is the canonical realization
of $A$ of Example \ref{bataninrealcan}.
Hence the functor $\nerf_\oper{A}\Pi^A_\infty S$
is a right Quillen functor. As any topological space
is fibrant, this implies that the functor $\nerf_\oper{A}\Pi^A_\infty S$
preserves weak equivalences. Our result then follows trivially from
Corollary \ref{controperquillenequ},
Theorem \ref{thmberger2}, Proposition \ref{Acellfirbant},
and from the equivalence of categories \eqref{berger3}.
\end{proof}

\begin{paragr}
Let $K$ be the universal contractible $\omega$-operad
(see Theorem \ref{existinitcontractoper}).
We defined in \ref{defweakcat} the weak $\omega$-categories
to be the $\oper{K}$-algebras. The $\omega$-operad
$K$ is also a contractible $\omega$-operad with system of
compositions as defined in \cite[Definitions 8.2 and 8.4]{Ba}.
Following M. Batanin, a \emph{weak $\omega$-groupoid} is a
weak $\omega$-category in which every cell is weakly invertible in
the precise sense given by \cite[9.5]{Ba}. For example,
for any topological space $X$, $\Pi^K_\infty(X)$
is a weak $\omega$-groupoid; see \cite[Theorem 9.1]{Ba}. We denote by
$\ogrp$ the full subcategory of $\alg_\oper{K}$
whose objects are the weak $\omega$-groupoids.
A morphism of weak $\omega$-groupoids is a
\emph{weak equivalence} if it is a weak equivalence
of $\oper{K}$-algebras. We will write $\ho(\ogrp)$
for the localization of the category of
weak $\omega$-groupoids by the weak equivalences.
The functor $\Pi^K_\infty$ of \eqref{deffondinftygroupoid}
thus induces by Theorem \ref{goodgroupoid1} a functor
$$\Pi_\infty : \ho(\Top)\To\ho(\ogrp) \ . $$
\end{paragr}

\begin{cor}\label{goodgroupoid2}
The functor $\Pi_\infty : \ho(\Top)\To\ho(\ogrp)$
is faithful and conservative.
\end{cor}

\begin{proof}
The inclusion functor from the category of
weak $\omega$-groupoids to the category of
weak $\omega$-categories induces a functor
$$i : \ho(\ogrp)\To\ho(\alg_\oper{K}) \ . $$
But it follows from Theorem \ref{goodgroupoid1}
that the functor $i\Pi_\infty=\Pi^K_\infty$ is an equivalence of categories.
This implies obviously our result.
\end{proof}

\phantom{\cite{Ba,bas,Be,BK,CS,Ler2,Si}
\cite{Ci1,Ci3,dug1,dug2,joyal,MoA1,Ta}
\cite{DHK,GZ,GJ2,Qu1,Lei,Leibook}
\cite{St1,St2,Gr1,BF,Hir}}
\bibliography{bibgr}
\bibliographystyle{amsalpha}
\end{document}